\documentclass[12pt,a4paper]{article}
\usepackage[english]{babel}
\usepackage[utf8x]{inputenc}
\usepackage{amsmath}
\usepackage[margin=1in]{geometry}
\usepackage{graphicx,cite}
\usepackage[section]{placeins}
\usepackage{float}
\usepackage{ amssymb }
\usepackage{hyperref}
\usepackage{amsthm}
\usepackage{amsfonts}
\usepackage{epsfig}
\usepackage{ textcomp }
\usepackage{mathtools}
\usepackage{authblk}
\usepackage{etoolbox}

\newcommand{\bb}{\begin{equation}}
\newcommand{\ee}{\end{equation}}

 \newtheorem{thm}{Theorem}
 \newtheorem{conj}[thm]{Conjecture}
 \newtheorem{prop}[thm]{Proposition}
 \newtheorem{lem}[thm]{Lemma}

\DeclareMathOperator{\diam}{diam}
\title{Minimizing the Number of Unions}
\date{(June 2023)}
\author{\v{Z}arko Ran\dj elovi\'c}
\affil{\small{ Centre for Mathematical Sciences
\\
Wilberforce Road\\
Cambridge CB3 0WB, U.K.}}
\affil{\texttt{zr233@cam.ac.uk}}

\begin{document}
\maketitle

\begin{abstract}
For a given number of $k$-sets, how should we choose them so as to minimize the union-closed family that they generate? Our main aim in this paper is to show that, if $\mathcal{A}$ is a family of $k$-sets of size $\binom{t}{k}$, and $t$ is sufficiently large, then the union-closed family generated by $\mathcal{A}$ has size at least that generated by the family of all $k$-sets from a $t$-set. This proves (for this size of family) a conjecture of Roberts. We also make some related conjectures, and give some other results, including a new proof of the result of Leck, Roberts and Simpson that exactly determines this minimum (for all sizes of the family) when $k=2$.

\end{abstract}

\section{Introduction} 

Let $\mathcal{A}=\{A_1,...,A_n\}$ be a family of $k$-sets where the $A_i$ are distinct. We are interested in minimizing the size of the union-closed family generated by $\mathcal{A}$, namely
$$\text{\text{\textless}}\mathcal{A}\text{\textgreater}=\{ A_{i_1}\cup A_{i_2}...\cup A_{i_s}| 1\le s\le n,1\le i_1,i_2,..,i_s\le n\}.$$
 It is natural to imagine that it is best to take the sets to be `as close together as possible', so for example all $k$-sets from a $t$-set if the size is $\binom{t}{k}$. Note that the union-closed family generated by this
family has size $|[t]^{(\geq k)}|$, where as usual $[t]$ denotes the set
$\{1,2,...,t\}$ and for any set $X$ we write $X^{(\geq k)}$ for the family of all
subsets of $X$ of size at least $k$. Our main aim is to prove this when $t$ is large. 
\begin{thm}
Let $\mathcal{A}$ be a family of $k$-sets of size $\binom{t}{k}$. Then for
$t$ sufficiently large we have $|\text{\text{\textless}}\mathcal{A}\text{\textgreater}|\geq |[t]^{(\geq k)}|$.
\end{thm}
This was conjectured by Roberts \cite{roberts}. In fact, he conjectured a result that should hold for all sizes of $\mathcal{A}$ : we discuss this in the final section of the paper. In principle the answer could have depended on the size of the ground set, say $N$, but in fact it does not. The first possible case would be $k=2$, and for this Leck, Roberts and Simpson showed an exact result.  Initial segments of colex are best (colex on $\mathbb{N}^{(k)}$ is the order in which $A<B$ if $\max(A\backslash B)<\max(B\backslash A)$ where $\mathbb{N}^{(k)}$ is the collection of all subsets of $\mathbb{N}$ of size $k$). 
\begin{thm}
(Leck, Roberts and Simpson \cite{robertseck}) Let $\mathcal{A}$ be a family of $2$-sets, and let $\mathcal{B}$ be the first
$|\mathcal{A}|$ $2$-sets in the colex order on $\mathbb{N}^{(2)}$. Then $|\text{\text{\textless}}\mathcal{A}\text{\textgreater}|\ge|\text{\text{\textless}}\mathcal{B}\text{\textgreater}|$. 
\end{thm}

 We give a new proof of this result. When $k\ge 3$ colex is not best. We discuss this at the end of the paper. We will use standard notation for set systems and graphs and their parameters. See Bollob\'{a}s \cite{10.5555/7228} for general background.

\section{Proofs of Results}

 We will start by giving an overview of the proof of Theorem $1$. Let $N$ be the size of the ground set. The idea of the proof is that we will first reduce the ground set to size close to $t$. To reduce the size of the ground set we first show that it cannot be "too big" and then we remove points belonging to few sets. If $N$ is not close to $t$ then after removing points belonging to few sets we get a set $X$ of size $N'$ none of whose elements belong to too few of the remaining sets. Then we show that most subsets of $X$ are indeed unions. We can show that $N'$ is close to $t$. Then unless $N$ is close to $t$ starting with $X$  we can add elements from the ground set one by one making many more unions. This will give more unions than the lower bound we are trying to prove. For the close to $t$ case,  similar to above removing elements belonging to few sets we get a set $Y$ such that nearly all subsets of $Y$ are unions, but this time we will have $|Y|\ge t$. We will get too many unions unless we have both $|Y|=t$ and $Y$ being the whole ground set which is then the $[t]^{(k)}$ case.  \\
 
 We will first prove a few lemmas. For the following lemmas and the proof of Theorem 1 we will always use the notation that $n=\binom{t}{k}$ where $n,t,k$ are integers such that $t\ge k>1$ and that $\mathcal{A}=\{ A_1,..,A_n \}$ is our family of $k$-sets. We will also define $G=\cup_{i=1}^nA_i$. We may view $G$ as the "ground set".\\
\\
 The following lemma will be used to show that the ground set size cannot be too big.
\begin{lem}
If $s\in \mathbb{N}$ and $|G|\ge sk$ then $|\text{\text{\textless}}\mathcal{A}\text{\textgreater}|\ge 2^s-1$.
\end{lem}
\textit{Proof.} We will prove a more general statement. Suppose that $l\in \mathbb{N}$ and $B_1,..,B_l$ are subsets of $\mathbb{N}$ such that $|B_i|\le k$ for all $i$. Let $T=\cup_{i=1}^lB_i$ and suppose $|T|\ge sk$. Finally let $$\mathcal{F}=\{ B_{i_1}\cup B_{i_2}...\cup B_{i_r}| 1\le r\le l,1\le i_1,i_2,..,i_r\le l\}.$$ We will show that $|\mathcal{F}|\ge 2^s-1$. We will prove that there is some $S\subset T$ such that $|S|=s$ and for any non-empty $X\subset S$ there is some $B\in \mathcal{F}$ such that $B\cap S=X$. This will prove the above statement and hence the lemma. To show this we will induct on $s$.\\
\\
For $s=1$ it is trivial as we may choose any $x\in B_1$ and set $S=\{x\}$.\\
\\
 Suppose it is true for some $s-1\ge 1$. Now consider the sets $B_1,..,B_l$ and $T$. For every $x\in T$ we define $T_x=\{i|x\in B_i\}$. Pick an $x\in T$ with the smallest $|T_x|$. Now without loss of generality $x\in B_1$. Let $T'=T\backslash B_1$. We have $|T'|\ge (s-1)k$. Now consider the collection of sets $$\mathcal{H}=\{B_i\cap T'|i\not \in T_x\}.$$ Suppose that $y\in T'$. We cannot have $T_x=T_y$ since $y\not \in B_1$ but $|T_x|\le |T_y|$ so there must be some $i\in T_y\backslash T_x$. But then $B_i\cap T'\in \mathcal{H}$ so $y\in \cup_{B\in \mathcal{H}}B $. Therefore $$\cup_{B\in \mathcal{H}}B=T'.$$
 
 Now by induction there is some $S'\subset T'$ of size $s-1$ such that for every non-empty $X\subset S'$ there is some $B$ which is a union of some sets in $\mathcal{H}$ such that $B\cap S'=X$. In other words the projection of a union-closed family containing $\mathcal{H}$ onto $S'$ contains all non-empty subsets of $S'$. From the definition of $\mathcal{H}$ this means that there is some $B\in \mathcal{F}$ such that $B\cap S'=X$ and $x\not \in B$. Now consider all possible unions of sets in the family $$\{B_i|i\not \in T_x\}\cup \{B_1\}$$ and let $S=S'\cup\{ x\}$. We get that the projection of $\mathcal{F}$ onto $S$ contains all non-empty subsets of $S$ and $|S|=s$. This completes the induction step and the proof. \hfill\qed
 \\
 \\
 We can see that $|[t]^{(\geq k)}|<2^t$ and so if $|\text{\text{\textless}}\mathcal{A}\text{\textgreater}|<|[t]^{(\geq k)}|$ then by Lemma 3 we must have $|G|\le tk$. This is good as we have restricted the ground set which will make it easier to work with. The next step is showing that we can essentially remove all elements in the ground set that are in too few sets. For any $x\in G$ define $d_x=|\{i|x\in A_i\}|$. 
 Also let $$\mathcal{A}_x=\{A\subset G|x\in A,|A|\ge k,\forall i\  (x\in A_i\Rightarrow A_i\not \subset A)\}.$$ 
 In other words $\mathcal{A}_x$ is the set of all subsets of $G$ of size at least $k$ for whom $x$ is one of the reasons that they are not in $\text{\textless}\mathcal{A}\text{\textgreater}$. \\
 \\
 The following lemma will be used to show that after removing elements in too few sets, most subsets of the remaining elements are in fact in $\text{\textless}\mathcal{A}\text{\textgreater}$.
 \begin{lem}
 Suppose that $s\in \mathbb{N}$ and $x\in G$ such that $d_x\ge s\binom{|G|}{k-2}$. Then we must have $|\mathcal{A}_x|\le 2^{|G|-s}.$
 \end{lem}
 \textit{Proof.} Consider the family $\mathcal{T}=\{B_i|x\in A_i, B_i=A_i\backslash\{x\}\}$. We have that $|\mathcal{T}|=d_x$. We may assume that $G\backslash \{ x\}=[p]$ for some $p$. Now we have $\mathcal{T}\subset [p]^{(k-1)}$. Notice that $\mathcal{A}_x$ consists of all subsets of $G$ containing $x$ that do not have a subset in $\mathcal{T}$. For $0\le r\le p-1$ we define the upper shadow of a family $\mathcal{F}\subset [p]^{(r)}$  to be $$\partial_+\mathcal{F}=\{A\cup \{i\}|A\in \mathcal{F},i\in [p]\backslash A\}.$$
 
 Also define $\partial_+^k\mathcal{F}=\partial_+ \partial_+ ... \partial_+ \mathcal{F}$ where $\partial_+$ is applied $k$ times ($\partial_+^0 \mathcal{F}=\mathcal{F}$). For $N,r\in \mathbb{N}$ define the lex order $[N]^{(r)}$ to be the order in which $A<B$ if $\min(A\backslash B)<\min(B\backslash A)$. Now let $\mathcal{I}$ be the initial segment of lex on $[p]^{(k-1)}$ of size $|\mathcal{T}|$. For $1\le r\le p$ we define the lower shadow of a family $\mathcal{F}\subset [p]^{(r)}$ to be $$\partial\mathcal{F}=\{A\backslash \{i\}|A\in \mathcal{F},i\in A\}.$$ Consider the map $g: [p] \rightarrow [p]$ given by $g(i)=p+1-i$. If $F: \mathcal{P}([p])\rightarrow \mathcal{P}([p])$ is given by $F(A)=g(A^c)$ then $F$ is a bijection that takes any initial segment of lex into a corresponding initial segment of colex. For any $\mathcal{B}\subset [p]^{(r)}$ we have $F(\partial_+(\mathcal{B}))=\partial F(\mathcal{B})$. Applying the Kruskal-Katona theorem (Kruskal \cite{kruskal1963number}, Katona \cite{katona1987theorem} and see Bollobás \cite{10.5555/7228}) we see that $|\partial F(\mathcal{T})|\ge |\partial F(\mathcal{I})|$ and hence $|\partial_+ \mathcal{T}|\ge|\partial_+ \mathcal{I}|$. 
 \\
 
 Note that the upper shadow of an initial segment of lex is an initial segment of lex. We can now show by induction that for any $1\le r\le p-k+1$ we have that $|\partial_+^r \mathcal{T}|\ge |\partial_+^r \mathcal{I}|$. We have by above that this is true for $r=1$. If it is true for $r<p-k+1$ then if $\mathcal{I'}$ is the initial segment of lex of length $|\partial_+^r\mathcal{T}|$ on $[p]^{(k-1+r)}$ then we know that $|\mathcal{I'}|\ge |\partial_+^r\mathcal{I}|$. Also since the upper shadow of an initial segment of lex is an initial segment of lex we have that  $\partial_+^r\mathcal{I}$ is an initial segment of lex on $[p]^{(k-1+r)}$ and hence $\partial_+^r\mathcal{I}\subset \mathcal{I}'$. Now applying the Kruskal-Katona theorem $$|\partial_+^{(r+1)}\mathcal{T|}\ge |\partial_+\mathcal{I'}|\ge |\partial_+^{(r+1)}\mathcal{I}|.$$ 
So by induction $|\partial_+^r \mathcal{T}|\ge |\partial_+^r \mathcal{I}|$ for all $1\le r\le p-k+1$. Now if we define the total upper shadow of a family $\mathcal{F}\subset [p]^{(r)}$ $$U_\mathcal{F}=\{A\subset [p]|X\subset A\  \text{for some}\ X\in \mathcal{F}\}$$
 then by above $|U_\mathcal{T}|\ge |U_\mathcal{I}|$. Notice also that $|\mathcal{A}_x|=|[p]^{(\ge k-1)}|-|U_\mathcal{T}|$. Since we want an upper bound on $|\mathcal{A}_x|$ we just need to show that $U_\mathcal{I}$ contains almost all sets in $[p]^{(\ge k-1)}$.  \\
 \\
 We know that in $[p]$ the number of $(k-1)$-sets containing a fixed element is $\binom{p-1}{k-2}$. We also know that $p=|G|-1$ So the number of $(k-1)$-sets containing at least one element from $[s]$ is at most $s\binom{|G|}{k-2}\le d_x$ which also means that $s\le p$. Note that in lex these sets are all before sets not containing any elements in $[s]$. so we must have that $\mathcal{I}$ contains all sets with at least one element from $[s]$ which means that the complement of $U_\mathcal{I}$ in $[p]^{(\ge k-1)}$ is a subset of $\mathcal{P}([p]\backslash[s])$. So we have $$|\mathcal{A}_x|\le |[p]^{(\ge k-1)}|-|U_\mathcal{I}|\le |\mathcal{P}([p]\backslash[s])|\le 2^{|G|-s}.$$ \hfill\qed
\\
\\
The following will be a useful fact about getting many new unions.
\begin{lem}
Suppose that $\mathcal{H}$ is a family of sets and $A$ is a set such that $A\not \subset \cup_{S\in \mathcal{H}}S$. Then $|\mathcal{H}\cup \{A\cup S|S\in \mathcal{H}\}|\ge(1+ \frac{1}{2^{|A|-1}})|\mathcal{H}|$.
\end{lem}
\textit{Proof.} Suppose that $x\in A\backslash \cup_{S\in \mathcal{H}}S$. If two sets $S_1,S_2$ in $\mathcal{H}$ have $A\cup S_1=A\cup S_2$ then $S_1$ and $S_2$ can only differ on $A\cap (\cup_{S\in \mathcal{H}}S)$ which is a fixed set of size at most $|A|-1$. So at most $2^{|A|-1}$ different sets in $\mathcal{H}$ can give the same union with $A$. This means that $$|\{A\cup S|S\in \mathcal{H}\}|\ge \frac{1}{2^{|A|-1}}|\mathcal{H}|.$$ However, no set containing $A$ is in $\mathcal{H}$ so we get the desired result.\hfill\qed
\\
\\
We now come to the heart of the proof.
\begin{lem}
For any real $D_k\ge 1$ there is a large enough constant $B$, depending only on $k$ and $D_k$, such that if $t>B$ and $|G|<t+D_k\log t$ then $|\text{\textless}\mathcal{A}\text{\textgreater}|\ge 2^t-\sum_{i=0}^{k-1}\binom{t}{i}$ and equality holds if and only if $\mathcal{A}=X^{(k)}$ where $X$ is a set with $|X|=t$.
\end{lem}
\textit{Proof.} Given $D_k$ suppose that $t>B$ and $|G|<t+D_k\log t$ where $B$ will be chosen later. Starting with $G$ and all of the $A_i$ we remove elements from $G$ one by one (and hence removing any of the $A_i$ containing those elements) if they are in less than $3\log(2tk)\binom{2t}{k-2}$ of the remaining sets $A_i$ until we get a $G'$ which satisfies that no $x\in G'$ is in less than $3\log(2tk)\binom{2t}{k-2}$ of the remaining sets. Notice that if $|G'|\le t-1$ then at some point in our process we were left with $G''\subset G$ such that $|G''|=t-1$. We can see that at that point we have removed at most $D_K\log t+1$ elements. This means that we have not removed that many sets either. In fact we have removed at most $3(D_k\log t+1)\log(2tk)\binom{2t}{k-2}$ sets. We will make sure that $B>3k$ to have that $\log(2tk)=\log t+\log(2k)\le 2\log t$, $D_k\log t+1\le 2D_k\log t$ and $2t\le 3(t-k)$. Now we just need to make sure that $$\frac{t-1}{\log^2 t}>12t^{1/2}\cdot 3^{(k-2)}(k-1)D_k$$ which is possible as $(t-1)/t^{1/2}\ge \sqrt{t}-1=e^{\frac{1}{2}\log t}-1> \log^3 t/48$ so we ensure that $B>e^{600\cdot 3^{(k-2)}(k-1)D_k}$. Now the number of sets removed when we are left with $G''$ is at most $$3(D_k\log t+1)\log(2tk)\binom{2t}{k-2}\le \frac{12\cdot 3^{(k-2)}D_k\log^2 t(t-k)^{k-2}}{(k-2)!}< t^{-1/2}\binom{t-1}{k-1}.$$ 

This means that for large enough $B$ we could have only removed less than $\binom{t-1}{k-1}$ sets and hence when we get $G''$ we still have more than $\binom{t}{k}-\binom{t-1}{k-1}=\binom{t-1}{k}$ sets remaining. This is a contradiction as $|G''|=t-1$ and hence we have that $|G'|\ge t$. Now we will show that if $B$ is large enough we can make sure that the vast majority of subsets of $G'$ are indeed in $\text{\textless}\mathcal{A}\text{\textgreater}$. Define for any $x\in G'$ $$\mathcal{A}_x'=\mathcal{A}_x\cap \mathcal{P}(G').$$ 
From the proof of Lemma 4, we can see that the lemma is true for any $n$, not just $n=\binom{t}{k}$. We will apply Lemma 4
on $G'$ and the remaining sets when we are left with $G'$. For large enough $B$ we have $|G'|\le |G|<2t$ so if we define $d_{G'}(x)=|\{A_i| x\in A_i, A_i\subset G'\}|$ we have $d_{G'}(x)\ge 3\log(2tk)\binom{|G'|}{k-2}$ for all $x\in G'$. Now by Lemma 4 $$|\mathcal{A}_x'|\le 2^{|G'|-[3\log (2tk)]}\le 2^{|G'|-2\log(2tk)}\le \frac{2^{|G'|}}{2tk}.$$
Now since every set in $G'^{(\ge k)}$ that is not in $\text{\textless}\mathcal{A}\text{\textgreater}$ must be in some $\mathcal{A}_x'$ we have by the union bound \begin{align}|\text{\textless}\mathcal{A}\text{\textgreater}|\ge |G'^{(\ge k)}|-|\cup_{x\in G'}\mathcal{A}_x'|\ge2^{|G'|}-\sum_{i=0}^{k-1}\binom{|G'|}{i}-|G'|\frac{2^{|G'|}}{2tk}.\end{align}
Since $|G'|\ge t>3k$ we have $\binom{|G'|}{i}<\frac{1}{2}\binom{|G'|}{i+1}$ for $0\le i\le k-1$ because $|G'|-i>2(i+1)$. Since exponential beats polynomial, for large enough $t$ we have $|G'|^k<\frac{2^{|G'|}}{2^{k+1}}$. We will take a large enough $B$ to have that $\sum_{i=0}^{k-1}\binom{|G'|}{i}\le \binom{|G'|}{k}<|G'|^k<\frac{2^{|G'|}}{2^{k+1}}$.
Now if $|G'|>t+1$ by $(1)$ we have $$|\text{\textless}\mathcal{A}\text{\textgreater}|\ge 2^{|G'|}-\frac{2^{|G'|}}{2^{k+1}}-\frac{2^{|G'|}}{2}\ge \frac{2^{|G'|}}{4}\ge 2^t$$ so we are done and if $|G'|=t+1$ then $3|G'|<4t\le 2tk$ so again by $(1)$ $$|\text{\textless}\mathcal{A}\text{\textgreater}|\ge 2^{|G'|}-\frac{2^{|G'|}}{2^{k+1}}-\frac{2^{|G'|}}{3}\ge \bigg(1-\frac{1}{8}-\frac{1}{3}\bigg)2^{|G'|}>\frac{1}{2}2^{|G'|}= 2^t$$ and we are also done. So we may assume that $|G'|=t$ and that $G'=[t]$. Similar to when we considered $G''$ we have now removed at most $\frac{1}{t^{1/2}}\binom{t-1}{k-1}$ sets and since $n=\binom{t}{k}$ we know that at most $\frac{1}{t^{1/2}}\binom{t-1}{k-1}$ sets in $[t]^{(k)}$ are not in $\mathcal{A}$ so for any $x\in [t]$ we have $$d_{G'}(x)\ge \frac{\sqrt{t}-1}{\sqrt{t}}\binom{t-1}{k-1}\ge \frac{t-1}{2k}\binom{t-2}{k-2}\ge \lceil t^{1/2}\rceil\binom{t-2}{k-2}$$ for large enough $B$. Now by Lemma 4 on $G'$ we have that $$|\mathcal{A}_x'|\le 2^{t-t^{1/2}}<\frac{2^t}{t2^{k+1}}$$ for large enough $B$ since $t^{1/2}>\frac{\log t}{\log 2} +(k+1)$ for large enough $t$. Now we have that $$|\cup_{x\in G'}\mathcal{A}_x'|\le \sum_{x\in G'}|\mathcal{A}_x'|\le \frac{2^t}{2^{k+1}}.$$ 
From before we have $\sum_{i=0}^{k-1}\binom{t}{i}\le \frac{2^t}{2^{k+1}}$ and thus $$|\text{\textless}\mathcal{A}\text{\textgreater}\cap [t]^{(\ge k)}|\ge 2^t-\frac{2^t}{2^{k+1}}-\frac{2^t}{2^{k+1}}=2^t\bigg(1-\frac{1}{2^k}\bigg).$$
Now if $G\neq G'$ there is some $A_i\in \mathcal{A},A_i\not \subset G'$ so by Lemma 5
$$|\text{\textless}\mathcal{A}\text{\textgreater}|\ge \bigg(1+\frac{1}{2^{k-1}}\bigg)|\text{\textless}\mathcal{A}\text{\textgreater}\cap [t]^{(\ge k)}|\ge \bigg(1+\frac{1}{2^{k-1}}\bigg)2^t\bigg(1-\frac{1}{2^k}\bigg)>2^t$$ and we are done. If $G=G'$ we must indeed have that $\mathcal{A}=[t]^{(k)}$ and so $|\text{\textless}\mathcal{A}\text{\textgreater}|=2^t-\sum_{i=0}^{k-1}\binom{t}{i}$. By the above we only have equality when $\mathcal{A}=X^{(k)}$ for some set $X$ where $|X|=t$. This proves the lemma. \hfill\qed
 \\
 
 We now prove our main result. All we need to do is reduce the ground set size enough to apply Lemma $6$. To do this we will apply Lemmas $3,4$ and $5$.
\\

\textit{Proof of Theorem 1.} We will show that there is a $C_k$ such that if $t>C_k$ we have $|\text{\textless}\mathcal{A}\text{\textgreater}|\ge |[t]^{(\geq k)}|=2^t-\sum_{i=0}^{k-1}\binom{t}{i}$. If $|G|\ge tk$ then by Lemma 3 we have $|\text{\textless}\mathcal{A}\text{\textgreater}|\ge 2^t-1$ so we may assume that $|G|<tk$. Now let $s=2\lceil\log(2tk)\rceil$. Keep removing elements from $G$ one by one (and all of the sets $A_i$ containing them) until no element that is left is in less than $s\binom{tk}{k-2}$ of the remaining sets. Suppose that $G'$ is the set of the remaining elements. We have removed at most $tks\binom{tk}{k-2}$ sets. Since $2k([t/2]+i)\ge tk$ for all $i\ge 1$ we have that $$tk\binom{tk}{k-2}\le \frac{(tk)^{k-1}}{(k-2)!}\le (k-1)(2k)^{k-1}\binom{[t/2]+k}{k-1}.$$

We will make sure $C_k>2k$. Now let $p=2\lceil\log(2tk)\rceil(k-1)(2k)^{k-1}<3\log t(2k)^{k}$. Then we have removed at most $p\binom{[t/2]+k}{k-1}$ sets. But we have also removed at least $\binom{t}{k}-\binom{|G'|}{k}$ sets. If $t> m\ge [t/2]+k$ then $$\binom{t}{k}-\binom{m}{k}=\sum_{l=m}^{t-1}\binom{l}{k-1}\ge (t-m)\binom{[t/2]+k}{k-1}.$$ Since exponential beats polynomial take $C_k$ large enough so that $$t-\bigg[\frac{t}{2}\bigg]-k\ge \frac{t}{3}=\frac{e^{\log t}}{3}>3\log t(2k)^{k}>p.$$
Then if $|G'|\le [t/2]+k$ then we have removed at least $$\binom{t}{k}-\binom{[t/2]+k}{k}\ge (t-[t/2]-k)\binom{[t/2]+k}{k-1}>p\binom{[t/2]+k}{k-1}$$ sets which is impossible. So we must have $|G'|\ge [t/2]+k$ and in fact if $|G'|=t-r$ then we must have $r\le p<3\log t(2k)^k$. We will first show that most subsets of $G'$ are actually in $\text{\textless}\mathcal{A}\text{\textgreater}$. \\
\\
Notice that if $x\in G'$ then $d_{G'}(x)\ge s\binom{tk}{k-2}\ge s\binom{|G'|}{k-2}$. This means we can apply Lemma 4 to $G'$. Define $\mathcal{A}_x'$ as in Lemma 6. We obtain that for all $x\in G'$ we have $$|\mathcal{A}_x'|\le 2^{|G'|-s}\le \frac{2^{|G'|}}{2^{2\log(2tk)}}\le \frac{2^{|G'|}}{2tk}.$$ 
This means that $$|\cup_{x\in G'}\mathcal{A}_x'|\le \sum_{x\in G'}|\mathcal{A}_x'|\le |G'|\frac{2^{|G'|}}{2tk}\le 2^{|G'|-1}.$$ 
If $A\subset G'$ and $|A|\ge k$ then $A\not \in \text{\textless}\mathcal{A}\text{\textgreater}$ if and only if there is some $x\in G'$ such that $A\in A_x'$. This means that 
\begin{align} |\text{\textless}\mathcal{A}\text{\textgreater}\cap \mathcal{P}(G')|=2^{|G'|}-|\cup_{x\in G'}\mathcal{A}_x'|-\sum_{i=0}^{k-1}\binom{|G'|}{i}\ge  2^{|G'|}-2^{|G'|-1}-\sum_{i=0}^{k-1}\binom{|G'|}{i}.\end{align}

We can easily bound the sum of the binomials since if $t>6k+2$ then we have  $|G'|>t/2>3k+1$ so just like in the proof of Lemma 6 we will take $C_k$ large enough so that $\sum_{i=0}^{k-1}\binom{|G'|}{i}\le \binom{|G'|}{k}<|G'|^k< 2^{|G'|-k-1}\le 2^{|G'|-2}$ since exponential beats polynomial. Now from (2) we have that $|\text{\textless}\mathcal{A}\text{\textgreater}\cap \mathcal{P}(G')|\ge 2^{|G'|-2}$. We will now show that there is a constant $D_k$ dependent on $k$ such that if $|G|\ge t+D_k\log t$ then we must have $|\text{\textless}\mathcal{A}\text{\textgreater}|\ge 2^t$. Let $\mathcal{H}_0=\text{\textless}\mathcal{A}\text{\textgreater}\cap \mathcal{P}(G')$ and $G_0=G'$. First of if $|G'|\ge t+2$ then we are done so suppose that $|G'|<t+2$. We will now show that we can get a lot more unions by adding sets with new elements. If we can pick some $x\in G\backslash G'$ and a set $A_i$ containing $x$ then notice that by taking  $\mathcal{H}_1=\mathcal{H}_0\cup \{A\cup A_i|A\in \mathcal{H}_0\}$ and $G_1=G_0\cup A_i$ we see that by Lemma 5  $$|\mathcal{H}_1|\ge \bigg(1+\frac{1}{2^{k-1}}\bigg)|\mathcal{H}_0|.$$

This is useful, because we have multiplied the total number of unions by a constant bigger than $1$ but dependent on $k$ and we have added at most $k$ new elements. We may keep on going as long as there are new elements and construct $\mathcal{H}_2,\mathcal{H}_3,...,\mathcal{H}_q$ and $G_2,G_3,...,G_q$ with $|H_i|\ge (1+\frac{1}{2^{k-1}})|H_{i-1}|$ and $|G_i\backslash G_{i-1}|\le k$ for $1\le i\le q$  and $G_q=G$. This means that if there are too many elements in $G$ we will get that the total number of unions is too big. First of notice that $|G'|=t-r>t-3\log t(2k)^k$. Now suppose that
$|G|\ge t+2+(2k+3k\log t(2k)^k)\log_{1+\frac{1}{2^{k-1}}}2$. This means that we can repeat the above process of adding a set and less than $k$ new elements at least $q\ge (2+3\log t(2k)^k)\log_{1+\frac{1}{2^{k-1}}}2$ times so we have that $$|\text{\textless}\mathcal{A}\text{\textgreater}|\ge |\mathcal{H}_0|\bigg(1+\frac{1}{2^{k-1}}\bigg)^{(2+3\log t(2k)^k)\log_{1+\frac{1}{2^{k-1}}}2}\ge 2^{|G'|-2}\cdot 2^{2+3\log t(2k)^k}\ge 2^t.$$
This means that there are constants $E_k$ and $D_k\ge 1$ dependent only on $k$ such that if $t>E_k$ and $|G|\ge t+D_k\log t$ we have that $|\text{\textless}\mathcal{A}\text{\textgreater}|\ge 2^t$. If $|G|<t+D_k\log t$ then since $D_k$ depends only on $k$ using Lemma 6 we get that there is a $B$ dependent only on $k$ such that if $t>B$ and $|G|< t+D_k\log t$ we have the desired result. Taking $C_k=\max (E_k,B)$ proves the theorem. \hfill\qed
\\

 In the rest of this section we will prove Theorem 2 which solves the case $k=2$ completely. We will show that if $n\in \mathbb{N}$ and $t\ge 2$ is the smallest positive integer such that $n\le \binom{t}{2}$ then $f(n,2)=2^t-2^{\binom{t}{2}-n}-t$. We note that $t=\lfloor \sqrt{2n}+\frac{3}{2}\rfloor$. 
\\

\textit{Proof of Theorem 2.} As above $n=|\mathcal{A}|$ and let $G$ be the graph with edge set $\mathcal{A}=\{A_1,..,A_n\}$ and vertex set $A_1\cup A_2...\cup A_n$. Suppose that $|\text{\textless}\mathcal{A}\text{\textgreater}|\le |\text{\textless}\mathcal{B}\text{\textgreater}|$. Let $t\in \mathbb{N}$ be such that $\binom{t-1}{2}<n\le \binom{t}{2} $. We see that $|G|\ge t$. Notice that $\mathcal{B}\subset [t]^{(2)}$ so $|\text{\textless}\mathcal{B}\text{\textgreater}|\le 2^t-t-1<2^t-1$. For any vertex $x\in G$ let $d(x)$ be the degree of $x$ in $G$. By considering all edges with $x$ we can make at least $2^{d(x)}-1$ distinct unions. Thus we know that the maximal degree is bounded by $\Delta (G)< t$. We also may assume that $\diam G\le 2$ since if there are some $x,y\in G$ such that the distance between them is $d(x,y)\ge 3$ then we may replace the family  $A_1,..,A_n$ with the family obtained from those sets by identifying $x$ and $y$. This will keep the size of the family to be $n$ and will not increase $|\text{\textless}\mathcal{A}\text{\textgreater}|$ as it will just identify $x,y$ in all of the unions. Denote by $\Gamma (x)$ the set containing $x$ and all elements adjacent to $x$ in $G$. We prove the following claim. \\
\\
\textbf{Claim.} For every $x\in G$ we have $|G|-|\Gamma (x)|<t$.\\
\\
\textit{Proof of Claim.} Suppose $x\in G$. Since $\diam G\le 2$ for every $y\in G \backslash \Gamma (x)$ there is a $z\in \Gamma (x)$ such that $yz\in E(G)$. This means if we consider the family $$\mathcal{C}=\{A\cap (G\backslash \Gamma(x))|A\in \text{\textless}\mathcal{A}\text{\textgreater}\}$$ then we have that each singleton subset of $G\backslash \Gamma (x)$ is in $\mathcal{C}$. But $\mathcal{C}$ is closed under unions and also has size at most $|\text{\textless}\mathcal{A}\text{\textgreater}|<2^t-1$ meaning that $|G\backslash \Gamma (x)|<t$ which proves the claim. \hfill\qed
\\
\\
With the bound on $\Delta (G)$ this gives $|G|\le 2t-1$. We will show that $|G|-t$ must be small. Just like in the proof of Theorem 1 instead of counting how many sets are in $\text{\textless}\mathcal{A}\text{\textgreater}$ we will instead count how many are in $\mathcal{F}=\mathcal{P}(G)\backslash \text{\textless}\mathcal{A}\text{\textgreater}$. Observe that if $S\in \mathcal{F}$ and $S\neq \emptyset$ then for some $x\in S$ for all $y\in S\backslash\{x\}$ we have $xy\not \in E(G)$. So for each set in $\mathcal{F}$ there is some element which prevents it from being in $\text{\textless}\mathcal{A}\text{\textgreater}$. Similar to our above definition of $\mathcal{A}_x$ let $$\mathcal{T}_x=\{S\subset G\ | x\in S,\forall y\in S\backslash\{x\}\ xy\not \in E(G)\}.$$
Denote by $s_x$ the degree of $x$ in $G^c$. We have $|\mathcal{T}_x|=2^{s_x}$ and by the above claim $s_x<t$ for any $x$. Note that $\mathcal{F}=(\cup_{x\in G}\mathcal{T}_x)\cup \{\emptyset\}$ so we have the bound \begin{align}|\mathcal{F}|\le 1+\sum\limits_{x\in G} |\mathcal{T}_x|\le1+ \sum\limits_{x\in G} 2^{s_x}.\end{align} 
We also know that $\sum\limits_{x\in G} s_x= \binom{|G|}{2}-n$. If $a\ge b>0$ then clearly $2^{a+1}+2^{b-1}>2^a+2^b$. Now consider variables $s_x'\in \mathbb{Z}_{\ge 0}$ for $x\in G$. If there are distinct $x,y\in G$ such that $0<s_x'\le s_y'<t-1$ then by above replacing $s_x',s_y'$ with $s_x'-1,s_y'+1$ increases $\sum\limits_{x\in G} 2^{s_x'}$. So the biggest the sum $\sum\limits_{x\in G} 2^{s_x'}$ could be subject to the constraints \begin{align}
\sum\limits_{x\in G} s_x'= \binom{|G|}{2}-n\end{align} \begin{align}s_x'<t\end{align} is if all except at most one of the $s_x'$ are either $0$ or $t-1$. Let $\binom{|G|}{2}-n=d(t-1)+m$ with non-negative integers $d,m$ and $0\le m<t-1$. Since the $s_x$ satisfy constraints (4) and (5) we must have $|G|\ge d$.
By (3) this means that \begin{align}|\mathcal{F}|\le d2^{t-1}+2^m+|G|-d.\end{align} Note that $|\text{\textless}\mathcal{A}\text{\textgreater}|=2^{|G|}-|\mathcal{F}|$ thus \begin{align}2^t-t-1\ge |\text{\textless}\mathcal{A}\text{\textgreater}|\ge 2^{|G|}-d2^{t-1}-2^m-|G|+d.\end{align}
We now want to bound $d$ in terms of $|G|$ and $t$. We know that \begin{align*} d\le \frac{\binom{|G|}{2}-n}{t-1}\le\frac{\binom{|G|}{2}-\binom{t-1}{2}}{t-1}=\frac{(|G|-t+1)(|G|+t-2)}{2(t-1)}\\ 
\\
\le \frac{(|G|-t+1)3(t-1)}{2(t-1)}=\frac{3}{2}(|G|-t+1).\end{align*}
We can also get a lower bound on $d$ since $$d=\Bigl\lfloor \frac{\binom{|G|}{2}-n}{t-1}\Bigr\rfloor\ge\Bigl\lfloor \frac{\binom{|G|}{2}-\binom{t}{2}}{t-1}\Bigr\rfloor=\Bigl\lfloor \frac{(|G|-t)(|G|+t-1)}{2(t-1)}\Bigr\rfloor\ge\Bigl\lfloor \frac{(|G|-t)2(t-1)}{2(t-1)}\Bigr\rfloor=|G|-t.$$
Now rearranging the terms in (7) and using the lower bound on $d$ we obtain $$(d+3)2^{t-1}\ge 2^t+2^m+d2^{t-1}\ge 2^{|G|}+1+t-|G|+d\ge 2^{|G|}.$$ Now let $l=|G|-t+1$. From the upper bound on $d$ and the previous inequality we obtain $$\frac{3}{2}l+3\ge d+3\ge 2^l$$ so we have $3l+6\ge 2^{l+1}$. Note that for $l=3$ this does not hold and by induction if $2^{l+1}>3l+6$ then $2^{l+2}>2(3l+6)>3(l+1)+6$ so we must have $l\le 2$.\\
\\
Now suppose that $l=2$. Then $|G|=t+1$. Since $\binom{t+1}{2}-\binom{t}{2}\le \binom{t+1}{2}-n<\binom{t+1}{2}-\binom{t-1}{2}$ we have that $t\le d(t-1)+m< 2t-1$ so either $d=1$ or $m=0,d=2$ which means that $d2^{t-1}+2^m\le 2^t+1$. If $d=1$ then $d2^{t-1}+2^m\le 2^t$ so in any case $d2^{t-1}+2^m-d\le 2^t-1$ Now from (6) we have $|\mathcal{F}|\le 2^t-1+(t+1)=2^t+t$. So we have $$|\text{\textless}\mathcal{A}\text{\textgreater}|=2^{t+1}-|\mathcal{F}|\ge 2^{t+1}-2^t-t=2^t-t>2^t-t-1\ge |\text{\textless}\mathcal{A}\text{\textgreater}|$$ which is a contradiction. \\
\\
Since $|G|\ge t$ we must have $l=1$. We may assume that $G=[t]$. We now have that $|\text{\textless}\mathcal{A}\text{\textgreater}|$ is minimized exactly when $|\mathcal{F}|$ is maximized. Now let $r=\binom{t}{2}-n$. We know that $\sum_{i=1}^ts_i=r$ but if we look at $\mathcal{T}_1,...,\mathcal{T}_t$ we see that at least $r$ sets appear in two of the $\mathcal{T}_i$. This is because all sets corresponding to the edges in $G^c$ appear in two of the $\mathcal{T}_i$. This means that \begin{align}|\mathcal{F}|=|(\cup_{x\in G}\mathcal{T}_x)\cup \{\emptyset\}|\le \sum_{i=1}^t2^{s_i}-r+1.\end{align}
We consider the sum $X=\sum_{i=1}^t2^{s_i}$ over all possible $G^c$ with $r$ edges. Without loss of generality we may assume that  $s_1\le s_2\le...\le s_t$. Now suppose that some edge in $G^c$ is not incident with $t$. Removing this edge we decrease $X$ by at most $2^{s_t}$ (since two of the $s_i$ decrease by 1. Now we can add an edge incident to $t$ (that is not already present) because $r< \binom{t}{2}-\binom{t-1}{2}=t-1$. This increases $X$ by at least $2^{s_t}+1$. So overall we have strictly increased $X$. This means that any $G^c$ which maximizes $X$ has all edges incident to the highest degree vertex so it is a star. In this case $X=2^r+r+t-1$.
So we have that $$|\mathcal{F}|\le 2^r+r+t-1-r+1=2^r+t.$$ Now we have that \begin{align*}|\text{\textless}\mathcal{A}\text{\textgreater}|\ge 2^t-|\mathcal{F}|\ge 2^t-2^r-t.\end{align*}
Notice that a set $S\subset [t]$ is not in $\text{\textless}\mathcal{B}\text{\textgreater}$ if and only if $|S|\le 1$ or $S=\{t\}\cup S'$ where $S'\subset \{t-r,t-r+1,...,t-1\}$. This means that $|\text{\textless}\mathcal{B}\text{\textgreater}|=2^t-t-1-(2^r-1)=2^t-2^r-t$. Thus $|\text{\textless}\mathcal{A}\text{\textgreater}|\ge |\text{\textless}\mathcal{B}\text{\textgreater}|$. This completes the proof. \hfill\qed 
\\
\\
We have now shown that $$f(n,2)=2^t-2^{\binom{t}{2}-n}-t.$$ Which families are extremal? If we have $|\text{\textless}\mathcal{A}\text{\textgreater}|=|\text{\textless}\mathcal{B}\text{\textgreater}|$ then by above we must have $|G|=t$ and we must maximize the sum in (8) so $G^c$ must be a star which means that this is the only case up to isomorphism when we have equality. Therefore the extremal families are only those isomorphic to $\mathcal{B}$.
\\
\\
\section{Related results}
A conjecture of Roberts \cite{roberts} asserts that Theorem 1 should be true for all $t$, not just for $t$ large.
\begin{conj}(Roberts \cite{roberts})
Let $k,t\in \mathbb{N}$ where $t\ge k$. Let $n=\binom{t}{k}$ and let $\mathcal{A}=\{A_1,..,A_n\}$ be a family of $n$ distinct $k$-sets. Then $|\text{\textless}\mathcal{A}\text{\textgreater}|\ge |[t]^{(\geq k)}|.$
\end{conj}
In a different direction, what happens for values between the binomial coefficients? Suppose that $\binom{t}{k}<n<\binom{t+1}{k}$. It is natural to assume that to minimize $|\text{\textless}\mathcal{A}\text{\textgreater}|$ over all families $\mathcal{A}$ of $n$ distinct $k$-sets we can pick a family on the ground set $[t+1]$ that contains $[t]^{(k)}$ and some sets containing $t+1$. We can see that our proof of Theorem 1 will give that if $|\text{\textless}\mathcal{A}\text{\textgreater}|$ is minimized then the size of the ground set $G$ must be close to $t$ if $t$ is sufficiently large. We can also see that for most values of $n$ we can apply the idea of Lemma 6 as well to show that the ground set must have size $t+1$. With these ideas in mind we are able to solve a few special cases when $n$ is very close to $\binom{t}{k}$. If $n=\binom{t}{k}-l$ where $l<\binom{k+1}{2}$ then for sufficiently large $t$ (even larger than we needed for $\binom{t}{k}$) following the same reasoning as before if $|\text{\textless}\mathcal{A}\text{\textgreater}|=f(n,k)$ then we may assume that $\mathcal{A}\subset [t]^{(k)}$. \\

Now we will show that $|\text{\textless}\mathcal{A}\text{\textgreater}|$ is minimized when $\mathcal{A}=[t]^{k)}\backslash \mathcal{C}_l'$ where  $\mathcal{C}_l$ is the initial segment of colex on $[t-1]^{(k-1)}$ of length $l$ and $\mathcal{C}_l'=\{A\cup \{t\}|A\in \mathcal{C}_l\}$ . First notice that the only sets in $[t]^{(\ge k)}$ that are not in $\text{\textless}\mathcal{A}\text{\textgreater}$ of size at least $k+1$ must have size exactly $k+1$. We then show the following proposition.
\begin{prop}
Let $l<\binom{k+1}{2}$ be a positive integer and $\mathcal{X}$ be a family of $k$-sets where $|\mathcal{X}|=l$. Let $b_{\mathcal{X}}$ be the number of $(k+1)$-sets with the property that at least $k$ of their $k$-subsets are in $\mathcal{X}$. Then $l\ge b_{\mathcal{X}}k-\binom{b_{\mathcal{X}}}{2}$.
\end{prop}
\textit{Proof.} Let $C_1,..,C_{b_{\mathcal{X}}}$ be the   $(k+1)$-sets counted in $b_{\mathcal{X}}$. Then if $i\neq j$ then $C_i$ and $C_j$ have at most one subset of size $k$ in common. Let $\mathcal{X}=\{X_1,..,X_l\}$. The number of pairs $(i,j)$ such that $X_i\subset C_j$ must be at least $b_{\mathcal{X}}k$ by counting for each $j$. On the other hand for any $j\neq j'$ there is at most one $i$ such that both $(i,j)$ and $(i,j')$ are counted. So $i$ is counted $k_i$ times where $\binom{k_i}{2}$ is the number of pairs $j\neq j'$ such that $C_j\cap C_{j'}=X_i$ and we have $\sum_{i=1}^l \binom{k_i}{2}\le \binom{b_{\mathcal{X}}}{2}$. Since $k_i\le 1+\binom{k_i}{2}$ by taking the sum over $i$ the total number of pairs counted is at most $l+\binom{b_{\mathcal{X}}}{2}$. This gives that $$l\ge kb_{\mathcal{X}}-\binom{b_{\mathcal{X}}}{2}$$ as required. \hfill\qed
\\

Going back to our problem notice that since $\mathcal{A}$ is $[t]^{(k)}$ with $l$ sets taken away we have that every set of size at least $k+1$ in $\mathcal{P}([t])$ is in $\text{\textless}\mathcal{A}\text{\textgreater}$ except those $(k+1)$-sets that have at least $k$ of its $k$-subsets removed. By the above proposition the number of those $(k+1)$-sets is at most the largest $s_l$ such that $l\ge s_lk-\binom{s_l}{2}=k+(k-1)+...+(k-s_l+1)$. Notice that when we have removed exactly the sets of $\mathcal{C}_l'$ from $[t]^{(k)}$ then that removes exactly $s_l$ sets of size $k+1$ from $\text{\textless}\mathcal{A}\text{\textgreater}$. This means that it is indeed optimal to take $[t]^{(k)}\backslash \mathcal{C}_l'$ for sufficiently large $t$.\\
\\

What about the case when $n=\binom{t}{k}+l$? In this case we can show that if $l$ is fixed then there is a large enough $t$ depending on both $k,l$ such that $|\text{\textless}\mathcal{A}\text{\textgreater}|$ is minimized when $$\mathcal{A}=[t]^{(k)}\cup \{Y_1,..,Y_l\}$$ where $Y_i=\{1,2,..,k-2,k-2+i,t+1\}$. Suppose that $|\text{\textless}\mathcal{A}\text{\textgreater}|=f(n,k)$. As before using Lemma 3 to get a bound on the ground set and then using the ideas in Lemma 6 and the proof of Theorem 1 we may assume that $\mathcal{A}$ contains $\binom{t}{k}-o(\binom{t-1}{k-1})$ sets in $[t]^{(k)}$ . Now suppose that $X_1,..,X_s$ are all of the sets in $\mathcal{A}$ that are not in $[t]^{(k)}$. Trivially $s\ge l$. Further we let $X_i'=X_i\cap [t]$ and $\mathcal{B}=\mathcal{A}\cap [t]^{(k)}$. Notice that if we have a set $B\in \text{\textless}\mathcal{B}\text{\textgreater}$ that contains one of the $X_i'$ then there is a set $X\in \text{\textless}\mathcal{A}\text{\textgreater}$ such that $X\not \subset [t]$ but $X\cap [t]=B$. This gives a lot of sets in $\text{\textless}\mathcal{A}\text{\textgreater}$ in addition to those that are in $\mathcal{P}([t])$.\\
\\
Let $\delta_s$ be the proportion of sets in $\mathcal{P}([t])$ that are in the total upper shadow of an initial segment of lex of length $s$ on $[t]^{(k-1)}$. Notice that (for $t$ sufficiently large) if $s=l$ then $\delta_s$ is actually the proportion of sets in $\mathcal{P}([(l+1)k])$ that are in the total upper shadow of that same initial segment of lex because that segment is inside $[(l+1)k]$. The same is true if $s=l+1$. Crucially $\delta_l$ and $\delta_{l+1}$ depend only on $k,l$ and notice that $\delta_{l+1}>\delta_l$. We see that if $\mathcal{C}=[t]^{(k)}\cup \{Y_1,..,Y_l\}$ we have that 
\begin{align}
    |\text{\textless}\mathcal{C}\text{\textgreater}|\sim 2^t(1+\delta_l)
    \end{align} for large $t$. Notice that the total upper shadow of $\mathcal{X}'=\{X_1',..,X_s'\}$ contains at least $l$ sets in $[t]^{(k-1)}$ where it is exactly $l$ if and only if $s=l$ and $|X_i'|=k-1$ for all $i$. If there are at least $l+1$ sets then applying the Kruskal-Katona theorem similar to the proof of Lemma 4 the total upper shadow of $\mathcal{X}'$ in $[t]$ has size at least $\delta_{t+1}2^t$. For large enough $t$ we can guarantee that $|\text{\textless}\mathcal{B}\text{\textgreater}|=(1-o(1))2^t$. This gives that there are at least $(\delta_{l+1}-o(1))2^t$ sets in $\text{\textless}\mathcal{A}\text{\textgreater}$ that are not subsets of $[t]$. Finally from (9) we have that $$|\text{\textless}\mathcal{A}\text{\textgreater}|\ge (1+\delta_{l+1}-o(1))2^t>|\text{\textless}\mathcal{C}\text{\textgreater}|\ge f(n,k)\ge |\text{\textless}\mathcal{A}\text{\textgreater}|$$ which is a contradiction. \\
    \\
    Thus we must have that $s=l$ and $|X_i'|=k-1$ for all $i$. Now let $a_i$ be the element in $X_i\backslash X_i'$. This means that $\mathcal{B}=[t]^{(k)}$. Now we know that each set in the total upper shadow of $\mathcal{X}'$ gives rise to at least $1$ set in $\text{\textless}\mathcal{A}\text{\textgreater}$ that is not contained in $[t]$. So \begin{align*}
    |\text{\textless}\mathcal{A}\text{\textgreater}|\ge |[t]^{(\ge k)}|+|\mathcal{U}_{\mathcal{X}'}|
    \end{align*}
    where $\mathcal{U}_{\mathcal{X'}}$ is the total upper shadow of $\mathcal{X}'$ in $[t]$. If $Y_i'=Y_i\cap [t]$ and $\mathcal{Y'}=\{Y_1',..,Y_l'\}$ we get that $|\mathcal{U}_{\mathcal{X}'}|\ge |\mathcal{U}_{\mathcal{Y'}}|$. This means that $|\text{\textless}\mathcal{A}\text{\textgreater}|\ge |[t]^{(\ge k)}|+|\mathcal{U}_{\mathcal{Y}'}|=|\text{\textless}[t]^{(k)}\cup \{Y_1,..,Y_l\}\text{\textgreater}|$ for sufficiently large $t$ which is what we wanted to show. Notice that also if some $a_i\neq a_j$ then $[t]$ gives rise to both $[t]\cup \{a_i\},[t]\cup \{a_j\}$ in $\text{\textless}\mathcal{A}\text{\textgreater}$ so we have $|\text{\textless}\mathcal{A}\text{\textgreater}|>|\text{\textless}[t]^{(k)}\cup \{Y_1,..,Y_l\}\text{\textgreater}|$ which is a contradiction. Thus we must have all $a_i$ to be the same. \\
    
    We have therefore found extremal families for some $n$ that are very close to $\binom{t}{k}$ but we need $t$ to be large enough. Still, we do not know for most values in between the binomial coefficients or for small values of $t$. 
    \\
    
    We will finish with a conjecture for the general result for all $n$. Colex is not the right order. As noted by Leck, Roberts and Simpson \cite{robertseck} if we consider $$\mathcal{A}=[4]^{(3)}\cup \{\{1,2,5\},\{1,3,5\},\{1,4,5\}\}$$ and let $\mathcal{B}$ be the initial segment of colex on $\mathbb{N}^{(3)}$ of length $7$. Then $|\text{\textless}\mathcal{A}\text{\textgreater}|=12$, but $|\text{\textless}\mathcal{B}\text{\textgreater}|=13$.  If $\binom{t}{k}<n<\binom{t+1}{k}$ then as mentioned before for most $n$ if $|\text{\textless}\mathcal{A}\text{\textgreater}|=f(n,k)$ then the ground set of $\mathcal{A}$ has size $t+1$ so we may assume that $\mathcal{A}\subset [t+1]^{(k)}$. If we further assume that $[t]^{(k)}\subset \mathcal{A}$ then the we have $|\text{\textless}\mathcal{A}\text{\textgreater}|=|[t]^{(\ge k)}|+|\mathcal{U}_{\mathcal{X}'}|$ where $\mathcal{X}'$ is as above and $|\mathcal{X}'|=n-\binom{t}{k}$. As before this is minimized when $\mathcal{X}'$ is the initial segment of lex on $[t]^{(k-1)}$. So something like "colex but lex inside a given maximal element" could be the correct order. Such a mixed ordering has occurred in other problems as well (see Engel and Leck \cite{engel} and also Duffus, Howard and Leader \cite{Duffus_2019}). We will define the $\textit{max-lex}$ order on $\mathbb{N}^{(k)}$ to be the order in which $A<B$ if either $\max A<\max B$ or else $\max A=\max B$ and $\min (A\Delta B)\in A$. The following lovely conjecture of Roberts \cite{roberts} includes Conjecture 7 -- we strongly believe it is true.

\begin{conj} (Roberts \cite{roberts})
Let $k\in \mathbb{N}$ and $\mathcal{A}\subset \mathbb{N}^{(k)}$. If $\mathcal{B}$ is the initial segment of max-lex on $\mathbb{N}^{(k)}$ of size $|\mathcal{A}|$ then $|\text{\textless}\mathcal{A}\text{\textgreater}|\ge |\text{\textless}\mathcal{B}\text{\textgreater}|.$
\end{conj}

\bibliographystyle{plain}
\bibliography{main}

\end{document}